%% file: RevisedORL_2018_539_Lindensjo2019MarARXIV.tex
\documentclass[10pt, a4paper]{article}
\usepackage{color,graphicx}
\usepackage{amsmath} 

\usepackage{endnotes} 

\usepackage{url}	 
\usepackage{amsfonts}

\usepackage[final]{changes}

\definechangesauthor[name={Kristoffer}, color=blue]{kri}

\input{tmacros.tex}

\begin{document}

\title{A regular equilibrium solves the extended HJB system}

\author{Kristoffer Lindensj\"o\footnote{Department of Mathematics, Stockholm University, SE-106 91 Stockholm, Sweden; kristoffer.lindensjo@math.su.se}}   

\date{March 15, 2019}

\maketitle 
 
\vspace{-22pt}
\begin{abstract} Control problems not admitting the dynamic programming principle are known as time-inconsistent. The game-theoretic approach is to interpret such problems as intrapersonal dynamic games and look for subgame perfect Nash equilibria. A fundamental result of time-inconsistent stochastic control is a verification theorem saying that solving the extended HJB system is a sufficient condition for equilibrium. We show that solving the extended HJB system is a necessary condition for equilibrium, under regularity assumptions. The controlled process is a general It\^o diffusion.
\end{abstract}

\vspace{2pt}
\noindent {\footnotesize \noindent {\bf AMS MSC2010:} 49J20; 49N90; 93E20; 35Q91; 49L99; 91B51; 91G80.
\vspace{1pt}

\noindent {\small \footnotesize  {\bf Keywords:} Dynamic inconsistency, Extended HJB system, Subgame perfect Nash equilibrium, Time-inconsistent preferences, Time-inconsistent stochastic control.}}

 \vspace{-10pt}
\section{Introduction} \vspace{-5pt} \label{intro}
Consider a controlled process $X^{\bf u}$ with initial data $(t,x)$ and the problem of choosing a control ${\bf u}$ that maximizes  
\beq  \label{valuefunc}
J(t,x,{\bf u}): = \mathbb{E}_{t,x}\bracket{F(x,X^{\bf u}_T)} +  G(x,\mathbb{E}_{t,x}\bracket{X^{\bf u}_T}),
\eeq
where $F$ and $G$ are deterministic functions and $t<T<\infty$ for a constant $T$. 
This problem is inconsistent in the sense that if a control ${\bf u}$ is optimal for the initial data $(t,x)$ then ${\bf u}$ is generally not optimal for other initial data $(s,y)$, which means that the dynamic programming principle cannot generally hold. This type of inconsistency is known as time-inconsistency.
    
\vspace{3pt} \noindent  The game-theoretic approach is to view problem \eqref{valuefunc} from the\deleted[id=kri,remark={}]{ perceptive} \added[id=kri,remark={}]{perspective}
of a person who controls the process $X^{\bf u}$  but whose preferences change when $(t,x)$ changes.\deleted[id=kri,remark={}]{The person is interpreted as comprising versions of herself, one version for each $(t,x)$. These versions of the person are then viewed as agents playing an intrapersonal sequential game regarding how to control $X^{\bf u}$.}\added[id=kri,remark={}]{ Specifically, the problem is viewed as a sequential non-cooperative intrapersonal game regarding how to control $X^{\bf u}$; where each $(t,x)$ corresponds to one player. See \cite[p. 549]{tomas-disc} for a more comprehensive interpretation along these lines.}  
 The approach is formalized by the definition of a subgame perfect Nash equilibrium, which is refinement of the notion of a Nash equilibrium for dynamic games, see Definition \ref{eq-def} below. The game-theoretic approach to time-inconsistency was first studied in a seminal paper by Strotz \cite{strotz} in which utility maximization under non-exponential discounting is studied. Selten \cite{selten1965spieltheoretische,selten1975reexamination} gave the first definition of subgame perfect Nash equilibrium, relying on the approach of Strotz. 
 
\vspace{3pt} \noindent Time-inconsistent problems were first studied in finance and economics. The time-inconsistency in this literature is typically due to the economic notions of 
endogenous habit formation, 
mean-variance utility and non-exponential discounting. 
These types of problems can be formulated and studied in the framework of the present paper. We formulate simple examples and give references in Section \ref{styl-exam}.

\vspace{3pt} \noindent The first general results on the game-theoretic approach to Markovian time-inconsistent stochastic control are 
due to Bj\"ork et al. who around 2010 defined the extended HJB system --- 
which is a system of simultaneously determined PDEs and an extension of the standard HJB equation ---
and proved a verification theorem in a general It\^o diffusion setting, see the recently published \cite{tomas-continpubl}. An analogous treatment of discrete time Markovian time-inconsistent stochastic control is presented in \cite{tomas-disc}. 
Early papers in mathematical finance to study the game-theoretic approach to time-inconsistent problems are \cite{basak2010dynamic,ekeland-arxiv,ekeland2010golden,Ekeland} where PDE methods for specific time-inconsistent problems ---  that are similar to the general method relying on the extended HJB system of \cite{tomas-continpubl,tomas-disc} --- are developed. Recent publications that use different versions of the extended HJB system to study time-inconsistent stochastic control problems include 
\cite
{bensoussan2014time,
tomas-mean, 
chen2014optimal,
he2013optimal,
kronborg2015inconsistent,
li2015time,
li2013optimal,
li2015timepartial,zeng2016robust}. 
In \cite{boualem}, the equilibrium of a time-inconsistent control problem is characterized by a stochastic maximum principle. Time-inconsistent stopping problems are studied in e.g. 
\cite{bayraktar2018,christensen2017finding,christensen2018time,huang2018time}\deleted[id=kri,remark={}]{\cite{miller2017nonlinear}}.   We refer to \cite{tomas-disc,christensen2017finding,christensen2018time,pedersen2016optimal,pedersen2017optimal} for short surveys of the literature on time-inconsistent stochastic control.

\vspace{3pt} \noindent Time-inconsistent problems can also be studied using the notion of dynamic optimality defined in \cite{pedersen2016optimal,pedersen2017optimal} and the pre-commitment approach. In the present setting pre-commitment corresponds to finding a control that maximizes \eqref{valuefunc} for a fixed $(t,x)$. For a definition of dynamic optimality and a comparison of the different approaches to time-inconsistency see \cite{christensen2017finding,pedersen2016optimal,pedersen2017optimal}.

\vspace{3pt} \noindent In Section \ref{problem} we formulate the time-inconsistent stochastic control problem corresponding to \eqref{valuefunc} in more detail and give the definition of equilibrium. In Section \ref{hjb-systemsec} we define the extended HJB system and prove the main result Theorem \ref{thm-equilibrium-solvesHJB} which says that solving the extended HJB system is a necessary condition for equilibrium, under regularity assumptions. 
To illustrate the main result we study a simple example in Section \ref{application}. Section \ref{gencase} contains a more general version of the main result.

 \vspace{-10pt}\subsection{Reasons for time-inconsistency} \vspace{-5pt}\label{styl-exam}
To give an idea of the type of time-inconsistent problems that are typically studied in finance and economics we here formulate three simple examples. We also give references to where problems of these types are studied. 
For further descriptions of endogenous habit formation, mean-variance utility and non-exponential discounting, and references, see e.g. \cite{tomas-disc,christensen2017finding}.

\vspace{3pt} \noindent  \textbf{Endogenous habit formation}: Problems of this type are studied in e.g.  \cite{abel1990asset,tomas-continpubl,christensen2017finding,detemple1992optimal,englezos2009utility,pollak1976habit}. 
As a simple example, consider an investor who controls the evolution of the wealth process $X^{\bf u}$ by dynamically adjusting the corresponding portfolio weights, see  \cite{karatzas1} for a standard model. 
Suppose the terminal time utility of the investor is 
$F(x,X^{\bf u}_T)$, where 
$F(x,\cdot)$ 
is a standard utility function for each fixed current wealth level $x$. In this case \eqref{valuefunc} becomes
\beqno
J(t,x,{\bf u}): = \mathbb{E}_{t,x}\bracket{(F(x,X^{\bf u}_T)}.
\eeqno
From an economic point of view this may be interpreted as the investor dynamically updating a habitual preference regarding the wealth level.

\vspace{3pt} \noindent  \textbf{Mean-variance utility}: Problems of this type are studied in e.g.  
\cite{basak2010dynamic,bayraktar2018,bensoussan2014time,bielecki2005continuous,tomas-disc,tomas-mean,christensen2018time,Czichowsky,he2013optimal,kronborg2015inconsistent,li2015time,li2013optimal,li2015timepartial,pedersen2016optimal,pedersen2017optimal,zeng2016robust}. As an example, consider the model above but an investor with mean-variance utility corresponding to
\beqno
J(t,x,{\bf u}) = \mathbb{E}_{t,x}\bracket{X^{\bf u}_T} -\frac{\gamma}{2}\mbox{Var}_{t,x}\bracket{X^{\bf u}_T},
\enskip \mbox{ where $\gamma>0$.}
\eeqno
The interpretation is that the investor wants a large expected wealth but is averse to risk measured by wealth variance. The parameter $\gamma$ corresponds to risk aversion.

\vspace{3pt} \noindent \textbf{Non-exponential discounting}: Problems of this type are studied in e.g. 
\cite{tomas-continpubl,
tomas-disc,
chen2014optimal,
ekeland-arxiv,
ekeland2010golden,
Ekeland,
huang2018time,
krusell,
weibull}.  As an example, consider the model above but an investor with a standard utility function $F$ and a deterministic discounting function $\varphi$ which cannot be rewritten as a standard exponential discounting function. Considering the time-space process, \eqref{valuefunc} becomes in this case
\beqno
J(t,x,{\bf u}) = \mathbb{E}_{t,x}\bracket{\varphi(T-t)F(X^{\bf u}_T)}, 
\eeqno
where $\varphi:[0,\infty)\rightarrow [0,1]$ is non-increasing with $\varphi(0)=1$.

 \vspace{-10pt}
\section{Problem formulation} \vspace{-5pt}\label{problem}

\noindent Consider a stochastic basis ${\pa{\Omega,{\cal F},P,\underline{\cal F}}}$ where $\underline{\cal F}$ is the augmented filtration generated by a $d$-dimensional Wiener process $W$. Consider a constant time horizon $T<\infty$ and an $n$-dimensional controlled  SDE
\beq \label{cont-SDE}
dX_s = \mu(s,X_s,{\bf u}(s,X_s))ds + \sigma(s,X_s,{\bf u}(s,X_s))dW_s,  \enskip X_t = x, \enskip t \leq s \leq T,
\eeq 
where
${\bf u}:[0,T]\times\mathbb{R}^n \rightarrow \mathbb{R}^k$, and 
$\mu:[0,T]\times\mathbb{R}^n\times\mathbb{R}^k \rightarrow \mathbb{R}^n$ and $\sigma:[0,T]\times\mathbb{R}^n\times\mathbb{R}^k\rightarrow M(n,d)$ are continous and satisfy standard global Lipschitz and linear growth conditions, see e.g.  \cite[sec 5.2]{Karatzas2}. $M(n,d)$ denotes the set of $n \times d$ matrices.

We also consider a mapping $U$ that restricts the set of values that controls ${\bf u}$ may take, see Definition \ref{admiss-ctrl}. Throughout the present paper we suppose $U$ and the functions $F$ and $G$ in \eqref{valuefunc} satisfy the following assumption.

\bass \label{rewassum} {\normalfont $F:\mathbb{R}^n\times\mathbb{R}^n \rightarrow \mathbb{R}$ is continuous and $G:\mathbb{R}^n\times\mathbb{R}^n \rightarrow \mathbb{R}$ satisfies $G\in C^2(\mathbb{R}^n\times\mathbb{R}^n)$. 
The control constraint mapping $U:[0,T]\times \mathbb{R}^n \rightarrow 2^{\mathbb{R}^k}$ is\deleted[id=kri,remark={}]{ continuous} 
\added[id=kri,remark={}]{such that for each $(t,x)\in [0,T)\times\mathbb{R}^n$ and each $u\in U(t,x)$ there exists a continuous control $\textbf{u}$ with $\textbf{u}(t,x) = u$}.
}
\eass
\added[id=kri,remark={}]{Note that constant control constraint mappings, which are used in most applications, trivially satisfy the condition in Assumption \ref{rewassum}.}

\bdf \label{admiss-ctrl} {\normalfont The set of admissible controls is denoted by ${\bf U}$. A control 
 ${\bf u}$ is said to be admissible if: ${\bf u}(t,x) \in U(t,x)$ for each $(t,x)\in [0,T]\times \mathbb{R}^n$, and for each $(t,x)\in[0,T)\times \mathbb{R}^n$ the SDE \cref{cont-SDE} has a unique strong solution $X^{\bf u}$ with the Markov property satisfying 
$\mathbb{E}_{t,x}\bracket{|F(x,X^{\bf u}_T)|}<\infty$ and 
$\mathbb{E}_{t,x}\bracket{||X^{\bf u}_T||}<\infty$.
}
\edf
\bdf  \label{auxfunc}{\normalfont  For any ${\bf u} \in {\bf U}$ the auxiliary functions $f_{\bf u}:[0,T]\times \mathbb{R}^n\times\mathbb{R}^n \rightarrow \mathbb{R}$ and $g_{\bf u}:[0,T]\times \mathbb{R}^n \rightarrow \mathbb{R}^n$ are defined by
\beqno
f_{\bf u}(t,x,y) = \mathbb{E}_{t,x}\bracket{F(y,X^{\bf u}_T)} \enskip
\textnormal{and} \enskip
g_{\bf u}(t,x) = \mathbb{E}_{t,x}\bracket{X^{\bf u}_T}.
\eeqno}
\edf
We are now ready to define the subgame perfect Nash equilibrium for the time-inconsistent stochastic control problem \cref{valuefunc}.  Definition \ref{eq-def} is in line with the equilibrium definition in e.g. \cite{tomas-continpubl,tomas-disc} to which we refer for a further motivation. 

\bdf [Equilibrium] \label{eq-def} \quad
{\normalfont \bei
\item Consider  a point $(t,x)\in[0,T) \times \mathbb{R}^n$, two controls  ${\bf u},{\bf \hat{u}} \in {\bf U}$ and a constant $h>0$.  
Let
\[
{\bf u}_h(s,y):=\begin{cases} 
	{\bf u}(s,y),        &\;  \mbox{ for } t \leq s <t+ h, \enskip y \in \mathbb{R}^n\\ 
	{\bf \hat{u}}(s,y),  &\;  \mbox{ for } t+ h \leq s \leq T,\enskip y \in \mathbb{R}^n.  
	\end{cases}
\]

\item \label{eq-control} The control ${\bf \hat{u}} \in {\bf U}$ is said to be an equilibrium control if, for any point $(t,x)\in[0,T) \times \mathbb{R}^n$ and any ${\bf u} \in {\bf U}$, it satisfies the equilibrium condition
\beq \label{eq-control-equation} 
\liminf_{h\searrow 0} \frac{J(t,x,{\bf \hat{u}})-J(t,x,{\bf u}_h)}{h} \geq 0.
\eeq

\item If ${\bf \hat{u}}$ is an equilibrium control then $V_{\bf \hat{u}}$ defined by $V_{\bf \hat{u}}(t,x)= J(t,x,{\bf \hat{u}})$ is said to be the corresponding  equilibrium value function and the quadruple $({\bf \hat{u}},V_{\bf \hat{u}},f_{\bf \hat{u}},g_{\bf \hat{u}})$ is said to be the corresponding equilibrium. 
\ei} 
\edf
The following definition will be used throughout the present paper.
\bdf \label{def-objects}\enskip {\normalfont
\bei
\item The differential operator ${\bf A}^{\bf u}$, corresponding to \cref{cont-SDE}, is defined by 
\beqno 
{\bf A}^{\bf u} = \frac{\partial}{\partial t} + \sum_{i=1}^{n}\mu_i(t,x,{\bf u}(t,x))\frac{\partial}{\partial x_i} + 
\frac{1}{2}\sum_{i,j=1}^{n}\sigma\sigma^T_{ij}(t,x,{\bf u}(t,x))\frac{\partial^2}{\partial x_ix_j}.
\eeqno
\added[id=kri,remark={}]{Moreover, for any constant $u \in \mathbb{R}^k$ we define}
\beqno 
\added[id=kri,remark={}]{{\bf A}^u = \frac{\partial}{\partial t} + \sum_{i=1}^{n}\mu_i(t,x,u)\frac{\partial}{\partial x_i} + 
\frac{1}{2}\sum_{i,j=1}^{n}\sigma\sigma^T_{ij}(t,x,u)\frac{\partial^2}{\partial x_ix_j}.}
\eeqno

\item Placing the third variable as a superscript for a function $f:[0,T] \times \mathbb{R}^n \times \mathbb{R}^n \rightarrow \mathbb{R}$, i.e. $f^y(t,x) = f(t,x,y)$, means that $y$ is to be taken as a constant. For example, $f^y\in C^{1,2}([0,T)\times\mathbb{R}^n)$ means that $f(t,x,y)$ is continuously differentiable with respect to $t$ and twice continuously differentiable with respect to $x$ for a fixed $y$,  
and ${\bf A}^{\bf u}f^y(t,x)$ involves only derivatives with respect to $t$ and $x$. 
Moreover, ${\bf A}^{\bf u}f(t,x,x)$ should be interpreted as ${\bf A}^{\bf u} \bar{f}(t,x)$ with $\bar{f}(t,x):= f(t,x,x)$.  

\item For a function $g:[0,T]\times \mathbb{R}^n \rightarrow \mathbb{R}^n$ we write $g(t,x)=(g_{1}(t,x),...,g_{n}(t,x))^T$ and let ${\bf A}^{{\bf u}}g(t,x):= ({\bf A}^{{\bf u}}g_1(t,x),...,{\bf A}^{{\bf u}}g_n(t,x))^T$.

\item The operator ${\bf H}^{\bf u}$ is defined by
\begin{align}
{\bf H}^{\bf u}g(t,x)= G_y(x,g(t,x)) {\bf A}^{{\bf u}}g(t,x), \enskip \mbox{ where } G_y(x,y):=\frac{\partial G}{\partial y}(x,y).
\label{Hoperator}
\end{align}
\added[id=kri,remark={}]{${\bf H}^u$ is defined analogously}.
\item 
\beq \label{diamondNOTA}
G \diamond g(t,x):= G(x,g(t,x)).
\eeq

\ei
}\edf
We will use the observation that \cref{valuefunc}, Definition \ref{auxfunc}, Definition \ref{eq-def} and \cref{diamondNOTA} imply,
\begin{align}
 V_{\bf \hat{u}}(t,x) &= J(t,x,{\bf \hat{u}}) \nonumber\\
&= f_{\bf \hat{u}}(t,x,x)+G\diamond g_{\bf \hat{u}}(t,x).\label{assumauxvaleq}
\end{align}

  \vspace{-20pt}
\section{The main result}  \vspace{-5pt}
\label{hjb-systemsec}
The extended HJB system is system a of simultaneously determined PDEs which we here define in line with \cite{tomas-continpubl}. 
Remark \ref{clar-rem} clarifies what constitutes a solution to the extended HJB system. 

\bdf [Extended HJB system] \label{HJB-system-def}{\normalfont For  $(t,x,y) \in [0,T) \times \mathbb{R}^n \times \mathbb{R}^n$,
\begin{align}
   \begin{split}
\quad\quad\quad\quad\quad\quad\quad\quad\quad\quad\quad\quad\quad\quad\quad\quad\quad\enskip
    {\bf A}^{{\bf \bar{u}}} f^y(t,x) &= 0, \label{kolmogorov1}\\
    f^y(T,x) &= F(y,x), \\
    \end{split}\\
   \begin{split}
    \quad\quad\quad\quad\quad\quad\quad\quad\quad\quad\quad\quad\quad\quad\quad\quad\quad\enskip\enskip
		{\bf A}^{{\bf \bar{u}}} g(t,x) &= 0, \label{kolmogorov2} \\
    g(T,x) &= x, \\
    \end{split}\\
		 \begin{split}
   \textnormal{sup}_{u \in U(t,x)} \{
   {\bf A}^uV(t,x)
 - {\bf A}^uf(t,x,x) + {\bf A}^uf^x(t,x) &  \label{HJB1A}\\
   - {\bf A}^u G\diamond g(t,x) + {\bf H}^ug(t,x)\} &=0,\\
    \end{split}\\
 \begin{split}
\quad\quad\quad\quad\quad\quad\quad\quad\quad\quad\quad\quad\quad\quad\quad\quad
\quad\quad\enskip
V(T,x) &=  F(x,x) + G(x,x)\label{HJB1B}\\
 \end{split}
\end{align}
where 
\begin{align} 
{\bf \bar{u}}(t,x) \in \textnormal{arg max}_{u \in U(t,x)} \{
   {\bf A}^uV(t,x)
 - {\bf A}^uf(t,x,x) + {\bf A}^uf^x(t,x)  \quad  \quad \quad \nonumber \\
 - {\bf A}^u G\diamond g(t,x) + {\bf H}^ug(t,x)\} \label{argmax}
			 \deleted[id=kri,remark={}]{\mbox{\sout{$= 0$}}}.
\end{align}}
\edf

\brem\label{clar-rem} {\normalfont For a fixed function ${\bf \bar{u}}$ equations \cref{kolmogorov1} and \cref{kolmogorov2} are   Kolmogorov backward equations. For fixed functions $f$ and $g$ equation \cref{HJB1A}--\eqref{HJB1B} is an HJB equation. The non-standard attribute of \cref{kolmogorov1}--\cref{HJB1B} is that ${\bf \bar{u}}$, $f$ and $g$ are not fixed in this way. Instead, \cref{kolmogorov1}--\cref{HJB1B} is a system simultaneously determined through \cref{argmax}.
 Let us describe what constitutes a solution: 
If four functions 
$V:[0,T]\times\mathbb{R}^n\rightarrow\mathbb{R},
f:[0,T]\times\mathbb{R}^n\times\mathbb{R}^n\rightarrow\mathbb{R}$,
$g:[0,T]\times\mathbb{R}^n\rightarrow\mathbb{R}^n$
and ${\bf \bar{u}}:[0,T]\times\mathbb{R}^n\rightarrow\mathbb{R}^k$, where ${\bf \bar{u}}(t,x) \in U(t,x)$ for each $(t,x)\in [0,T]\times\mathbb{R}^n$,  
satisfy the following conditions then $({\bf \bar{u}}, V, f, g)$ is a solution to the extended HJB system:
\bei 
\item  $f^y$ and ${\bf \bar{u}}$   satisfy \cref{kolmogorov1}, for each fixed $y\in \mathbb{R}^n$.
\item  $g$ and ${\bf \bar{u}}$   satisfy \cref{kolmogorov2}.
\item $V$ satisfies \eqref{HJB1B}.
\item $V,f,g$ and ${\bf \bar{u}}$ satisfy \eqref{HJB1A} and \eqref{argmax}, i.e for each fixed $(t,x)\in[0,T)\times\mathbb{R}^n$ the inequality $
   {\bf A}^uV(t,x)
 - {\bf A}^uf(t,x,x) + {\bf A}^uf^x(t,x) 
 - {\bf A}^u G\diamond g(t,x) + {\bf H}^ug(t,x)\leq  0 $ holds for each constant $u\in U(t,x)$, and it holds with equality for the constant $u:={\bf \bar{u}}(t,x)$.
\ei}
\erem

\noindent In order to prove the main result, Theorem \ref{thm-equilibrium-solvesHJB}, we need Lemma \ref{feyn-kac-lemma}, Lemma \ref{feyn-kac-lemma2} and Proposition \ref{lemmanewtt} below. We remark that Lemma \ref{feyn-kac-lemma} and Lemma \ref{feyn-kac-lemma2} are versions of the Feynman-Kac formula. A proof is included for the sake of completeness. 
We will use the following definition. 

\bdf \label{space-L2} { \normalfont Consider a control  ${\bf u} \in {\bf U}$. 
For a function  
\deleted[id=kri,remark={}]{\emph{h}}$\added[id=kri,remark={}]{k}:[0,T]\times\mathbb{R}^n\rightarrow\mathbb{R}$  
we write 
\deleted[id=kri,remark={}]{\emph{h}}$\added[id=kri,remark={}]{k}\in L^2_T(X^{\bf u})$ 
if, for each $(t,x) \in [0,T)\times \mathbb{R}^n$, 
\added[id=kri,remark={}]{there exists a constant $\bar{h}>0$ satisfying $t+\bar{h}<T$ such that}
\begin{small}
\begin{align*}
\added[id=kri,remark={}]{
{
\mathbb{E}_{t,x}\bracket{
\sup_{0\leq h \leq \bar{h} }
\left|\int_t^{t+h} \frac{{\bf A}^{{\bf u}}k(s,X^{\bf u}_s)}{h}ds\right|
+\int_t^{t+\bar{h}} \left|\left|
\frac{\partial k}{\partial x} (s,X^{\bf u}_s) \sigma(s,X^{\bf u}_s, {\bf {u}}(s,X^{\bf u}_s))
\right|\right|^2ds
}<\infty.
}}\end{align*} 
\end{small}
 \begin{small}
\vspace{-3mm}
\deleted[id=kri,remark={}]{$\mathbb{E}_{t,x}\bracket{\int_t^T 
\left(
|{\bf A}^{{\bf u}}h(s,X^{\bf u}_s)|
+\left|\left|
\frac{\partial h}{\partial x} (s,X^{\bf u}_s) \sigma(s,X^{\bf u}_s, {\bf {u}}(s,X^{\bf u}_s))
\right|\right|
^2 
\right) ds
}<\infty.$} 
\end{small}
}\edf
\blem \label{feyn-kac-lemma} {\normalfont Consider a continuous control ${\bf u}\in {\bf U}$. Suppose the auxiliary function $f_{\bf {u}}$ satisfies $f_{\bf {u}}^y \in C^{1,2}([0,T) \times  \mathbb{R}^n)\cap L^2_T(X^{\bf u})$, for any fixed $y \in \mathbb{R}^n$. Then, $f_{\bf {u}}^y$ is, for any fixed $y \in \mathbb{R}^n$, a  solution to the PDE
\begin{align*}
    {\bf A}^{{\bf u}} f^y(t,x) = 0, \enskip 
    f^y(T,x) = F(y,x), 		\enskip (t,x) \in [0,T) \times  \mathbb{R}^n.
\end{align*}}
\elem
\noindent \emph{Proof.} By definition $f_{\bf u }^y(t,x) = \mathbb{E}_{t,x}[F(y,X_T^{\bf u})]$ and the boundary condition is therefore satisfied. Consider an arbitrary point $(t,x,y) \in [0,T) \times \mathbb{R}^n \times \mathbb{R}^n$. Let $X^{\bf {u}}$ be the strong solution to the SDE \cref{cont-SDE} for the initial data $(t,x)$. 
Consider an arbitrary constant\deleted[id=kri,remark={}]{ $h>0$, satisfying $t+h<T$}\added[id=kri,remark={}]{ $h$ with $0<h<\bar{h}$, where $\bar{h}$ is as in Definition \ref{space-L2}}. 
The Markov property and It\^{o}'s formula imply that
\begin{align*}
0 &=\mathbb{E}_{t,x}\bracket{f_{\bf {u}}^y({t+h},X^{\bf {u}}_{t+h})} -f_{\bf {u}}^y(t,x)\nonumber   \\
&=\mathbb{E}_{t,x}\bracket{
 \int_t^{t+h} {\bf A}^{\bf {u}}f_{\bf u}^y(s,X^{\bf {u}}_s)ds}, 
\end{align*}
where the It\^{o} integral vanished since $f_{\bf u}^y\in L^2_T(X^{\bf u})$.
Hence, 
\beq \label{FKeq}
\mathbb{E}_{t,x}\bracket{\frac{\int_t^{t+h} {\bf A}^{\bf {u}}f_{\bf {u}}^y(s,X^{\bf {u}}_s)ds}{h}} = 0.
\eeq
The condition $f^y_{\bf u}\in L^2_T(X^{\bf u})$ implies that we can use dominated convergence when sending sending $h\searrow 0$ in \cref{FKeq}. Moreover, the integrand in \cref{FKeq} is continuous in $s$ for a.e. $\omega$, since $\mu$, $\sigma$ and ${\bf u}$ are continuous and $f^y_{\bf u}\in C^{1,2}([0,T) \times  \mathbb{R}^n)$. Hence, 
$
\lim_{h \searrow 0}  \mathbb{E}_{t,x}\bracket{\frac{\int_t^{t+h} {\bf A}^{\bf u}f_{\bf {u}}^y(s,X^{\bf {u}}_s)ds}{h}}
= {\bf A}^{\bf {u}}f_{\bf {u}}^y(t,x).
$ 
The result follows. 
\hfill $\square$

\blem \label{feyn-kac-lemma2} {\normalfont Consider a continuous control ${\bf u}\in {\bf U}$. Suppose the elements of the auxiliary function $g_{{\bf u}}$ satisfy  $g_{{\bf u},i} \in C^{1,2}([0,T) \times  \mathbb{R}^n)\cap L^2_T(X^{\bf u})$, $i=1,...,n$.  Then $g_{\bf {u}}$ is a solution to the PDE
\begin{align*}
   \begin{split}
    {\bf A}^{{\bf u}} g(t,x) = 0, \enskip
    g(T,x) = x, \enskip \mbox{ for }(t,x) \in [0,T)  \times  \mathbb{R}^n. \\
    \end{split}
\end{align*}
}
\elem
\noindent \emph{Proof.} The proof is analogous to that  of Lemma \ref{feyn-kac-lemma} and is  omitted.
\hfill $\square$

\bprop \label{lemmanewtt} {\normalfont  Consider two controls ${\bf v},{\bf \tilde v} \in {\bf U}$ where ${\bf v}$ is continuous.  Suppose the auxiliary functions $f_{\bf {\tilde v}}$ and $g_{{\bf \tilde v}}$ satisfy  
$f_{\bf {\tilde v}}^y,g_{{\bf \tilde v},i} \in C^{1,2}([0,T) \times  \mathbb{R}^n)\cap L^2_T(X^{\bf v})$, 
for any fixed $y \in \mathbb{R}^n$ and  $i=1,...,n$.
Consider a point $(t,x)\in [0,T)\times\mathbb{R}^n$. Let
\[
{\bf v}_h(s,y):=\begin{cases} 
	{\bf v}(s,y),        &\;  \mbox{ for } t \leq s <t+ h, \enskip y \in \mathbb{R}^n\\ 
	{\bf \tilde v}(s,y),  &\;  \mbox{ for } t+ h \leq s \leq T,\enskip y \in \mathbb{R}^n.  
	\end{cases}
\]
Let $v:={\bf v}(t,x)$. Then,
\begin{align}  \label{newlemma11} 
\lim_{h \searrow 0} \frac{f_{\bf \tilde v}(t,x,x) -f_{{\bf v}_h}(t,x,x)}{h} & = -{\bf A}^vf_{\bf \tilde v}^x(t,x), \\
\label{newlemma22}
\lim_{h \searrow 0} \frac{G \diamond g_{\bf \tilde v}(t,x)-G \diamond g_{{\bf v}_h}(t,x)}{h} & = - {\bf H}^vg_{\bf \tilde v}(t,x).
\end{align}}
\eprop
\noindent \emph{Proof.} 
$\mathbb{E}_{t,x}\bracket{f_{\bf \tilde v}^x(t+h,X^{{\bf v}}_{t+h})} -f^x_{\bf\tilde  v}(t,x)
= \mathbb{E}_{t,x}\bracket{\int_t^{t+h} {\bf A}^{\bf v}f_{\bf \tilde v}^x(s,X^{{\bf v}}_s)ds}$ is found as in Lemma \ref{feyn-kac-lemma}.
By definition, ${\bf v}_h$ and ${\bf v}$ coincide on $[t,t+h]$, except at the point $t+h$.
By definition, ${\bf v}_h(s,y)$ and ${\bf \tilde  v}(s,y)$ coincide on $[t+h,T]$. Thus,
\begin{align*}
\mathbb{E}_{t,x}\bracket{f_{\bf \tilde v}^x({t+h},X^{{\bf v}}_{t+h})}
&= \mathbb{E}_{t,x}\bracket{\mathbb{E}_{t+h,X_{t+h}^{{\bf v}}}[F(x,X_T^{\bf \tilde v})]} \nonumber \\
&= \mathbb{E}_{t,x}\bracket{\mathbb{E}_{t+h,X_{t+h}^{{\bf v}_h}}[F(x,X_T^{{\bf v}_h})]}\nonumber \\
&= f_{{\bf v}_h}(t,x,x). 
\end{align*}
From the above \added[id=kri,remark={}]{it} follows \added[id=kri,remark={}]{that} 
$
f_{{\bf   v}_h}(t,x,x) - f^{x}_{\bf \tilde v}(t,x)= 
\mathbb{E}_{t,x}\bracket{\int_t^{t+h} {\bf A}^{{\bf v}}f_{\bf \tilde v}^x(s,X^{\bf v}_s)ds}.
$ 
Using arguments analogous to those in the proof of Lemma \ref{feyn-kac-lemma} we thus obtain
\begin{align*} 
\lim_{h \searrow 0} \frac{f_{\bf \tilde v}(t,x,x) -f_{{\bf v}_h}(t,x,x)}{h}
&= \lim_{h \searrow 0} \frac{ - \mathbb{E}_{t,x}\bracket{\int_t^{t+h} {\bf A}^{{\bf v}}f_{\bf \tilde  v}^x(s,X^{\bf v}_s)ds}
}{h}\\
&= -{\bf A}^{\bf v}f_{\bf \tilde v}^x(t,x),
\end{align*}
which, since $v:={\bf v}(t,x)$, means that \cref{newlemma11} holds. Using the same arguments as above we obtain
$
g_{{\bf v}_h}(t,x) = g_{\bf \tilde v}(t,x)+\mathbb{E}_{t,x}\bracket{\int_t^{t+h}{\bf A}^{{\bf v}}g_{\bf \tilde v}(s,X^{\bf v}_s) ds}.
$ 
Standard Taylor expansion gives
\begin{align*}
&G\left(x,g_{\bf \tilde v}(t,x)+\mathbb{E}_{t,x}\bracket{\int_t^{t+h}{\bf A}^{{\bf v}}g_{\bf \tilde v}(s,X^{\bf v}_s) ds} \right)\\
&\enskip= G\left(x,g_{\bf \tilde v}(t,x) \right)+
G_y\left(x,g_{\bf \tilde v}(t,x)\right)\mathbb{E}_{t,x}\bracket{\int_t^{t+h}{\bf A}^{{\bf v}}g_{\bf \tilde v}(s,X^{\bf v}_s) ds} +o(h).
\end{align*}
Hence, \cref{newlemma22} follows from,
\begin{align*}
&\lim_{h \searrow 0}\frac{G(x,g_{\bf \tilde v}(t,x))-G(x,g_{{\bf v}_h}(t,x))}{h}\\
&\enskip =\lim_{h \searrow 0}\frac{-G_y\left(x,g_{\bf \tilde v}(t,x)\right)\mathbb{E}_{t,x}\bracket{\int_t^{t+h}{\bf A}^{{\bf v}}g_{\bf \tilde v}(s,X^{\bf v}_s) ds} +o(h) }{h}\\
&\enskip= -G_y(x,g_{\bf \tilde v}(t,x)){\bf A}^{{\bf v}}
g_{\bf \tilde v}(t,x)
.
\end{align*} 
\hfill $\square$


\vspace{3pt} \noindent Let us now define what is meant by a regular equilibrium and present main result Theorem \ref{thm-equilibrium-solvesHJB}. An example with a regular equilibrium is studied in Section \ref{application}.
\bdf \label{regularEQ} {\normalfont An equilibrium $({\bf \hat{u}}, V_{\bf \hat{u}}, f_{\bf \hat{u}}, g_{\bf \hat{u}})$  is said to be regular if:
\renewcommand{\theenumi}{(\roman{enumi})}%
\begin{enumerate}

\item \label{defitem1} The equilibrium control ${\bf \hat{u}}$ is continuous.

\item \label{defitem2} $f_{\bf \hat{u}}^y, g_{{\bf \hat{u}},i} \in L^2_T(X^{\bf \hat u})$ and $f_{\bf \hat{u}}^y, g_{{\bf \hat{u}},i},\bar{f}\in C^{1,2}([0,T) \times  \mathbb{R}^n)$ 
for each fixed $y \in \mathbb{R}^n$ and $i=1,...,n$, where $\bar{f}(t,x):= f_{\bf \hat{u}}(t,x,x)$.

\item \label{defitem3} For each $(t,x) \in [0,T) \times  \mathbb{R}^n$ and each $u \in U(t,x)$, there exists a continuous control ${{\bf u}}\in{{\bf U}}$ with ${{\bf u}}(t,x)=u$ such that  $f_{\bf \hat{u}}^y, g_{{\bf \hat{u}},i} \in L^2_T(X^{\bf u})$.
\end{enumerate}
}
\edf
\deleted[id=kri,remark={}]{Regarding the technical condition \ref{defitem3} we remark that Assumption \ref{rewassum} ensures that for each $(t,x)$ and each $u\in U(t,x)$, there exists a continuous control ${{\bf u}}$ with ${\bf u} (t,x) = u$, and in particular this holds when the constraint mapping $U(t,x)$ is constant.}

\bth \label{thm-equilibrium-solvesHJB} {\normalfont A regular equilibrium $({\bf \hat{u}},V_{\bf \hat{u}},f_{\bf \hat{u}},g_{\bf \hat{u}})$ solves the extended HJB system.}
\eth
\noindent \emph{Proof.}  Lemma \ref{feyn-kac-lemma} implies that the auxiliary function 
$f^y_{\bf \hat{u}}(t,x)$ and the equilibrium control ${{\bf \hat{u}}}$  satisfy \cref{kolmogorov1}, for each $y\in \mathbb{R}^n$. 
Lemma \ref{feyn-kac-lemma2} implies that the auxiliary function $g_{\bf \hat{u}}(t,x)$ and ${{\bf \hat{u}}}$  satisfy \cref{kolmogorov2}.  Sufficient regularity for the use of these lemmas is provided by \ref{defitem1} and \ref{defitem2} in Definition \ref{regularEQ}. The boundary condition \cref{HJB1B} is directly verified using \cref{diamondNOTA}, \cref{assumauxvaleq} and Definition \ref{auxfunc}. 
Now consider an arbitrary point $(t,x)\in [0,T)\times\mathbb{R}^n$. In order to show that the equilibrium $({\bf \hat{u}}, V_{\bf \hat{u}} , f_{\bf \hat{u}} ,g_{\bf \hat{u}} )$ is a solution to the extended HJB system we only have left to show that the inequality \eqref{concmainpart343} below holds for any $u\in U(t,x)$ and that it holds with equality for $u:={\bf \hat{u}}(t,x)$:
\begin{align}
& {\bf A}^{u}V_{\bf \hat{u}}(t,x) - {\bf A}^{u}f_{\bf \hat{u}}(t,x,x) + {\bf A}^{u}f^x_{\bf \hat{u}}(t,x)  \nonumber \\
& \quad\quad\quad\quad\quad- {\bf A}^{u}G \diamond g_{\bf \hat{u}}(t,x) + {\bf H}^{u}g_{\bf \hat{u}}(t,x)
\leq 0.\label{concmainpart343}
\end{align}
Consider an arbitrary $u\in U(t,x)$. From \cref{assumauxvaleq} \added[id=kri,remark={}]{it} follows \added[id=kri,remark={}]{that} 
\begin{align}
{\bf A}^uV_{\bf \hat{u}}(t,x) 
= {\bf A}^uf_{\bf \hat{u}}(t,x,x)
+ {\bf A}^uG\diamond g_{\bf \hat{u}}(t,x),\label{asdasa44}
\end{align}  
where differentiability is provided by \ref{defitem2} and Assumption \ref{rewassum}. Consider a continuous control $\bf u$ satisfying ${\bf u}(t,x)=u$ for which $f_{\bf \hat{u}}^y, g_{{\bf \hat{u}},i} \in L^2_T(X^{\bf u})$, cf. \ref{defitem3}. Use Proposition \ref{lemmanewtt} and \cref{asdasa44} to find
\begin{align}
&\lim_{\added[id=kri,remark={}]{h}\deleted[id=kri,remark={}]{\mbox{\sout{$n$}}} \searrow 0} \frac{
f_{\bf \hat{u}}(t,x,x) + G \diamond g_{\bf \hat{u}}(t,x)
-(f_{{\bf u}_h}(t,x,x) + G \diamond g_{{\bf u}_h}(t,x))}{h}\quad\quad\quad\quad\quad\nonumber \\
&\enskip= - {\bf H}^ug_{\bf \hat{u}}(t,x) -{\bf A}^uf_{\bf \hat{u}}^x(t,x)\nonumber \\
&\enskip= -({\bf A}^{u}V_{\bf \hat{u}}(t,x) - {\bf A}^{u}f_{\bf \hat{u}}(t,x,x) + {\bf A}^{u}f^x_{\bf \hat{u}}(t,x) \nonumber \\
&\quad \quad \quad \quad \quad \quad \quad \quad   -{\bf A}^{u}G \diamond g_{\bf \hat{u}}(t,x) + {\bf H}^{u}g_{\bf \hat{u}}( t,x)). \label{asaw34r34}
\end{align}
Now use the definition of $J(t,x,{{\bf u}})$ 
and the assumption that ${\bf \hat{u}}$ is an equilibrium control, cf. the equilibrium condition \cref{eq-control-equation}, 
to obtain
\begin{align}
&\lim_{\added[id=kri,remark={}]{h}\deleted[id=kri,remark={}]{\mbox{\sout{$n$}}} \searrow 0} \frac{f_{\bf \hat{u}}(t,x,x) + G \diamond g_{\bf \hat{u}}(t,x)
-(f_{{\bf u}_h}(t,x,x) + G \diamond g_{{\bf u}_h}(t,x))}{h}\nonumber \\
&\enskip= \lim_{h \searrow 0} \frac{J(t,x,{\bf \hat{u}}) - J(t,x,{\bf u}_h)}{h} \geq 0.\label{asaw34r342}
\end{align}
Recall that $u\in U(t,x)$ was arbitrarily chosen. Hence, \eqref{asaw34r34} and \eqref{asaw34r342} imply that \cref{concmainpart343} holds for any $u\in U(t,x)$.

\vspace{3pt} \noindent Since $f_{\bf \hat{u}}^y$ and ${{\bf \hat{u}}}$ satisfy \cref{kolmogorov1} for any $y$ 
\added[id=kri,remark={}]{it} follows \added[id=kri,remark={}]{that} 
${\bf A}^{\bf \hat{u}}f^x_{\bf \hat{u}}(t,x)=0$. 
Since $g_{\bf \hat{u}}$ and ${{\bf \hat{u}}}$ satisfy \cref{kolmogorov2} 
\added[id=kri,remark={}]{it} follows \added[id=kri,remark={}]{that}  
$
{\bf A}^{\bf \hat{u}}g_{\bf \hat{u}}(t,x)=0
$ 
which with \cref{Hoperator} gives 
${\bf H}^{{\bf \hat{u}}}g_{\bf \hat{u}}(t,x) = 0$. 
From \cref{assumauxvaleq} \added[id=kri,remark={}]{it} follows  \added[id=kri,remark={}]{that} 
${\bf A}^{{\bf \hat{u}}}V_{\bf \hat{u}}(t,x)
= {\bf A}^{{\bf \hat{u}}}f_{\bf \hat{u}}(t,x,x) +  {\bf A}^{{\bf \hat{u}}}G \diamond g_{\bf \hat{u}}(t,x).
$ 
Hence,
\beqno \label{concmainpart2}
{\bf A}^{{\bf \hat{u}}}V_{\bf \hat{u}}(t,x) 
- {\bf A}^{{\bf \hat{u}}}f_{\bf \hat{u}}(t,x,x) + {\bf A}^{{\bf \hat{u}}}f^x_{\bf \hat{u}}(t,x)
- {\bf A}^{{\bf \hat{u}}}G \diamond g_{\bf \hat{u}}(t,x) + {\bf H}^{{\bf \hat{u}}}g_{\bf \hat{u}}(t,x)
=0.
\eeqno
This is equivalent to \cref{concmainpart343} holding with equality for $u:={\bf \hat{u}}(t,x)$.
\hfill $\square$

 \vspace{-10pt}
\subsection{An example}  \vspace{-5pt}\label{application}
Let us study a simple time-inconsistent problem to illustrate Theorem \ref{thm-equilibrium-solvesHJB}. 
Suppose a person controls the evolution of a one-dimensional diffusion process with constant volatility by choosing its drift function. Specifically,
\beqno \label{appdyn}
dX_t = {\bf u}(t,X_t)dt + \sigma dW_t,
\eeqno  
where admissible controls are restricted to the interval $U = [-a,a]$ for some $a>0,\sigma>0$. Suppose the person would like a large difference between the current value of the process and its value at a fixed terminal time $T$. Specifically,
\begin{align}
J(t,x,{\bf u})  = \mathbb{E}_{t,x}\bracket{(X^{\bf u}_T-x)^2}.\label{examplevaluefunc}
\end{align}
This corresponds to $F(x,y)=(y-x)^2$ and  $G(x,y) = 0$. We make the ansatz that ${\bf \hat{u}}= 0$ is an equilibrium control. Simple calculations give us the corresponding auxiliary functions
$g_{\bf \hat{u}}(t,x) = x$  and $f_{\bf \hat{u}}(t,x,y) = (x-y)^2 + \sigma^2(T-t)$.
Hence, 
$
\frac{\partial   f^y_{\bf \hat{u}}(t,x)}{\partial t} = -\sigma^2, 
\frac{\partial   f^y_{\bf \hat{u}}(t,x)}{\partial x} = 2x-2y$ and$ 
\frac{\partial^2 f^y_{\bf \hat{u}}(t,x)}{\partial x^2}  = 2.
$ 
Let us now show that ${\bf \hat{u}}=0$ does indeed satisfy the equilibrium condition \cref{eq-control-equation}. Consider an arbitrary control ${\bf u} \in {\bf U}$ and an arbitrary point $(t,x)$. 
From It\^{o}'s formula \added[id=kri,remark={}]{it} follows \added[id=kri,remark={}]{that} 
\begin{align*}
& \mathbb{E}_{t,x}\bracket{f_{\bf \hat{u}}^x({t+h},X^{{\bf u}}_{t+h})}-f^{x}_{\bf \hat{u}}(t,x) \\
&\enskip = \mathbb{E}_{t,x}\bracket{
\int_t^{t+h} {\bf A}^{{\bf u}} f_{\bf \hat{u}}^x(s,X^{{\bf u}}_s)ds+
\sigma \int_t^{t+h}   \frac{\partial f^x_{\bf \hat{u}}(s,X^{\bf u}_s)}{\partial x}   dW_s}\\
&\enskip=   \mathbb{E}_{t,x}\bracket{
\int_t^{t+h}{\bf u}(s,X^{\bf u}_s)(2X_s^{\bf u}-2x)ds}.
\end{align*}
Using arguments similar to those in the proof of Proposition \ref{lemmanewtt} we find
$
f_{{\bf u}_h}(t,x,x)-f_{\bf \hat{u}}(t,x,x) = 
\mathbb{E}_{t,x}\bracket{\int_t^{t+h}   {\bf u}(s,X^{\bf u}_s)(2X_s^{\bf u}-2x)ds}.
$ 
Since $G(x,y)=0$  \added[id=kri,remark={}]{it} follows \added[id=kri,remark={}]{that}  
\begin{align*}
J(t,x,{\bf \hat{u}}) - J(t,x,{\bf u}_h) 
&= f_{\bf \hat{u}}(t,x,x)  - f_{{\bf u}_h}(t,x,x)\nonumber\\
&= \mathbb{E}_{t,x}\bracket{\int_t^{t+h}{\bf u}(s,X^{\bf u}_s)(2x-2X_s^{\bf u})ds}\nonumber\\
&\geq -2a\mathbb{E}_{t,x}\bracket{\int_t^{t+h}|X_s^{\bf u}-x|ds}.
\end{align*}
It follows that the equilibrium condition  \cref{eq-control-equation} holds and that ${\bf \hat{u}}=0$ therefore is an equilibrium control. The corresponding equilibrium is 
\begin{align} \label{eqex22}
({\bf \hat{u}} ,V_{\bf \hat{u}} ,f_{\bf \hat{u}},g_{\bf \hat{u}} ) =(0,\sigma^2(T-t),(x-y)^2 + \sigma^2 (T-t),x).
\end{align}
From Theorem \ref{thm-equilibrium-solvesHJB} \added[id=kri,remark={}]{it} follows that \eqref{eqex22} solves the extended HJB system corresponding to 
$F(x,y)=(y-x)^2$ and  $G(x,y) = 0$, which is also easily verified.

\begin{remark} {\normalfont The pre-commitment approach is in this example to maximize \eqref{examplevaluefunc} over admissible controls, by treating $x$ as an arbitrary but fixed parameter. It is easy to see that the pre-commitment optimal control is, for any fixed $x$,
\[
{\bf u}(t,y)=\begin{cases} 
	a,        &\;  \mbox{ for } y \geq x \\ 
	-a,  &\;  \mbox{ for } y < x. 
	\end{cases}
\]}
\end{remark}

 \vspace{-10pt}\subsection{A more general problem} \vspace{-5pt} \label{gencase}
In this section we include a running time function $H$ and allow $F$ and $G$ to depend on the initial time $t$. Specifically, we consider
\begin{align*}
J(t,x,{\bf u})     
&:=  \mathbb{E}_{t,x}\bracket{
   \int_t^T H(t,x,r,X^{\bf u}_r,{\bf u}(r,X^{\bf u}_r))dr 
   +F(t,x,X^{\bf u}_T)}&\\ 
&\enskip\enskip\enskip +   G(t,x,\mathbb{E}_{t,x}\bracket{X^{\bf u}_T})
\end{align*}
where $F:[0,T]\times\mathbb{R}^n\times\mathbb{R}^n \rightarrow \mathbb{R}$ and $G:[0,T]\times\mathbb{R}^n\times\mathbb{R}^n \rightarrow \mathbb{R}$ satisfy conditions analogous to those in Assumption \ref{rewassum} and 
$H:[0,T]\times\mathbb{R}^n \times [0,T]\times \mathbb{R}^n \times \mathbb{R}^k \rightarrow \mathbb{R}$ is continuous and bounded. The definition of an admissible control is analogous to Definition \ref{admiss-ctrl}. Let
\begin{align*}
G\diamond g(t,x) &:= G(t,x,g(t,x)),\\
{\bf H}^{\bf u}g(t,x) &:= G_y(t,x,g(t,x)) {\bf A}^{{\bf u}}g(t,x),\\
f_{\bf u}(t,x,s,y) &:=  \mathbb{E}_{t,x}\bracket{
   \int_t^T H(s,y,r,X^{\bf u}_r,{\bf u}(r,X^{\bf u}_r))dr 
   +F(s,y,X^{\bf u}_T)}.
\end{align*}
The equilibrium definition is analogous to Definition \ref{eq-def}. Placing the third and fourth variables of a function $f:[0,T]\times\mathbb{R}^n\times[0,T]\times\mathbb{R}^n\rightarrow \mathbb{R}$ as superscripts, i.e. $f^{s,y}(t,x)= f(t,x,s,y)$, means $s$ and $y$ are to be taken as constant. 
\bdf [Extended HJB system II] \label{HJB-system-def2} {\normalfont For  $(t,x,s,y) \in [0,T) \times \mathbb{R}^n\times [0,T) \times \mathbb{R}^n$,
\begin{align*}
  \quad\quad\quad\quad\quad\quad\quad\quad\enskip
    {\bf A}^{{\bf \bar{u}}} f^{s,y}(t,x) + H(s,y,t,x,{\bf \bar{u}}(t,x)) &= 0, \\
   f^{s,y}(T,x) &= F(s,y,x), \\
      \quad\quad\quad\quad\quad\quad\quad\quad\quad\quad\quad\quad\quad\quad\quad\quad\quad\quad\enskip
		{\bf A}^{{\bf \bar{u}}} g(t,x) \enskip &= 0, \\
    g(T,x) &= x, \\
   \textnormal{sup}_{u \in U(t,x)} \{
   {\bf A}^uV(t,x)
 - {\bf A}^uf(t,x,t,x) + {\bf A}^uf^{t,x}(t,x) &   \\
   - {\bf A}^u G\diamond g(t,x) + {\bf H}^ug(t,x) + H(t,x,t,x,u) \} &=0,\\
  V(T,x) &=  F(T,x,x) + G(T,x,x),
  \end{align*}
where \begin{align*} 
{\bf \bar{u}}(t,x) \in \textnormal{arg max}_{u \in U(t,x)} \{
   {\bf A}^uV(t,x)
 - {\bf A}^uf(t,x,t,x) + {\bf A}^uf^{t,x}(t,x) \quad \quad \quad \quad \quad  \\
 - {\bf A}^u G\diamond g(t,x) + {\bf H}^ug(t,x)  +  H(t,x,t,x,u)\}\deleted[id=kri,remark={}]{\mbox{\sout{$= 0$}}}.
\end{align*}} 
\edf
The definition of a regular equilibrium is analogous to that of Definition \ref{regularEQ}. 
Theorem \ref{thm-equilibrium-solvesHJB2} generalizes the main result of this paper to the present setting. The proof is analogous to that of Theorem \ref{thm-equilibrium-solvesHJB} and is omitted.

\bth \label{thm-equilibrium-solvesHJB2} {\normalfont A regular equilibrium $({\bf \hat{u}},V_{\bf \hat{u}} ,f_{\bf \hat{u}} ,g_{{\bf \hat{u}}})$ solves the extended HJB system II.}
\eth
 
\noindent {\footnotesize \textbf{Acknowledgments:} The author is grateful to Tomas Bj\"ork and Jan-Erik Bj\"ork for helpful discussions. An early pre-print version of this paper had the title \emph{Time-inconsistent stochastic control: solving the extended HJB system is a necessary condition for regular equilibria}.}

 \vspace{-15pt} 

\bibliography{kristofferBibl}

\end{document}

%% file: tmacros.tex
\newcommand{\beq}{\begin{equation}}
  \newcommand{\eeq}{\end{equation}}
  \newcommand{\beqno}{\begin{displaymath}}
  \newcommand{\eeqno}{\end{displaymath}}
  \newcommand{\beqar}{\begin{eqnarray}}
  \newcommand{\eeqar}{\end{eqnarray}}
  \newcommand{\beqarno}{\begin{eqnarray*}}
  \newcommand{\eeqarno}{\end{eqnarray*}}

\newcommand{\cref}[1]{(\ref{#1})}

\newcommand{\ba}{\begin{array}}
\newcommand{\ea}{\end{array}}

\newcommand{\ben}{\begin{enumerate}}
  \newcommand{\en}{ \end{enumerate}}

\newcommand{\bei}{\begin{itemize}}
  \newcommand{\ei}{ \end{itemize}}

\newcommand{\bed}{\begin{description}}
  \newcommand{\ed}{\end{description}}

\newcommand{\bec}{\begin{center}}
  \newcommand{\ec}{\end{center}}

\newcommand{\bprop}{\begin{proposition}}
  \newcommand{\eprop}{\end{proposition}}

\newcommand{\bdf}{\begin{definition}}
  \newcommand{\edf}{\end{definition}}

\newcommand{\bth}{\begin{theorem}}
  \newcommand{\eth}{\end{theorem}}

\newcommand{\bcon}{\begin{con}}
  \newcommand{\econ}{\end{con}}

\newcommand{\bcor}{\begin{corollary}}
  \newcommand{\ecor}{\end{corollary}}

\newcommand{\bpr}{\begin{problem}}
  \newcommand{\epr}{\end{problem}}

\newcommand{\blem}{\begin{lemma}}
  \newcommand{\elem}{\end{lemma}}

\newcommand{\brem}{\begin{remark}}
  \newcommand{\erem}{\end{remark}}
 
\newcommand{\bass}{\begin{assumption}}
  \newcommand{\eass}{\end{assumption}}

\newcommand{\bres}{\begin{result}}
  \newcommand{\eres}{\end{result}}

\newcommand{\bexm}{\begin{example}}
  \newcommand{\eexm}{\end{example}}
  
  \newcommand{\bema}{\begin{main}}
  \newcommand{\ema}{\end{main}}

\newtheorem{theorem}{Theorem}[section]
  \newtheorem{assumption}[theorem]{Assumption}
  \newtheorem{definition}[theorem]{Definition}
  \newtheorem{corollary}[theorem]{Corollary}
  \newtheorem{lemma}[theorem]{Lemma}
  \newtheorem{proposition}[theorem]{Proposition}
 \newtheorem{result}[theorem]{Result}
  \newtheorem{example}[theorem]{Example}
  \newtheorem{"definition"}[theorem]{"Definition"}
  \newtheorem{remark}[theorem]{Remark}
  \newtheorem{con}[theorem]{Conjecture}
  \newtheorem{problem}[theorem]{Problem}
 \newtheorem{main}[theorem]{Main question}

\newcommand{\bracket}[1]{\left[ {#1} \right]}

\newcommand{\pa}[1]{\left( {#1} \right)}

\newcommand{\model}{(\Omega,{\cal F},\Pi,X,\underline{\cal F})}






%% file: RevisedORL_2018_539_Lindensjo2019MarARXIV.bbl
\begin{thebibliography}{10}

\bibitem{abel1990asset}
A.~B. Abel.
\newblock Asset prices under habit formation and catching up with the
  {J}oneses.
\newblock Technical report, National Bureau of Economic Research, 1990.

\bibitem{basak2010dynamic}
S.~Basak and G.~Chabakauri.
\newblock Dynamic mean-variance asset allocation.
\newblock {\em The Review of Financial Studies}, 23(8):2970--3016, 2010.

\bibitem{bayraktar2018}
E.~Bayraktar, J.~Zhang, and Z.~Zhou.
\newblock Time consistent stopping for the mean-standard deviation problem ---
  the discrete time case.
\newblock arXiv preprint 1802.08358, 2018.

\bibitem{bensoussan2014time}
A.~Bensoussan, K.~Wong, S.~C.~P. Yam, and S.-P. Yung.
\newblock Time-consistent portfolio selection under short-selling prohibition:
  From discrete to continuous setting.
\newblock {\em SIAM Journal on Financial Mathematics}, 5(1):153--190, 2014.

\bibitem{bielecki2005continuous}
T.~R. Bielecki, H.~Jin, S.~R. Pliska, and X.~Y. Zhou.
\newblock Continuous-time mean-variance portfolio selection with bankruptcy
  prohibition.
\newblock {\em Mathematical Finance}, 15(2):213--244, 2005.

\bibitem{tomas-continpubl}
T.~Bj{\"o}rk, M.~Khapko, and A.~Murgoci.
\newblock On time-inconsistent stochastic control in continuous time.
\newblock {\em Finance and Stochastics}, 21(2):331--360, 2017.

\bibitem{tomas-disc}
T.~Bj{\"o}rk and A.~Murgoci.
\newblock A theory of {M}arkovian time-inconsistent stochastic control in
  discrete time.
\newblock {\em Finance and Stochastics}, 18(3):545--592, 2014.

\bibitem{tomas-mean}
T.~Bj{\"o}rk, A.~Murgoci, and X.~Y. Zhou.
\newblock Mean-variance portfolio optimization with state-dependent risk
  aversion.
\newblock {\em Mathematical Finance}, 24(1):1467--9965, 2014.

\bibitem{chen2014optimal}
S.~Chen, Z.~Li, and Y.~Zeng.
\newblock Optimal dividend strategies with time-inconsistent preferences.
\newblock {\em Journal of Economic Dynamics and Control}, 46:150--172, 2014.

\bibitem{christensen2017finding}
S.~Christensen and K.~Lindensj{\"o}.
\newblock On finding equilibrium stopping times for time-inconsistent
  {M}arkovian problems.
\newblock {\em SIAM Journal on Control and Optimization}, 56(6):4228--4255,
  2018.

\bibitem{christensen2018time}
S.~Christensen and K.~Lindensj{\"o}.
\newblock On time-inconsistent stopping problems and mixed strategy stopping
  times.
\newblock {\em arXiv preprint arXiv:1804.07018}, 2018.

\bibitem{Czichowsky}
C.~Czichowsky.
\newblock Time-consistent mean-variance portfolio selection in discrete and
  continuous time.
\newblock {\em Finance and Stochastics}, 17(2):227--271, 2013.

\bibitem{detemple1992optimal}
J.~B. Detemple and F.~Zapatero.
\newblock Optimal consumption-portfolio policies with habit formation.
\newblock {\em Mathematical Finance}, 2(4):251--274, 1992.

\bibitem{boualem}
B.~Djehiche and M.~Huang.
\newblock A characterization of sub-game perfect equilibria for {SDE}s of
  mean-field type.
\newblock {\em Dynamic Games and Applications}, 6(1):1--27, 2015.

\bibitem{ekeland-arxiv}
I.~Ekeland and A.~Lazrak.
\newblock Being serious about non-commitment: subgame perfect equilibrium in
  continuous time.
\newblock arXiv preprint math/0604264, 2006.

\bibitem{ekeland2010golden}
I.~Ekeland and A.~Lazrak.
\newblock The golden rule when preferences are time inconsistent.
\newblock {\em Mathematics and Financial Economics}, 4(1):29--55, 2010.

\bibitem{Ekeland}
I.~Ekeland and T.~Pirvu.
\newblock Investment and consumption without commitment.
\newblock {\em Mathematics and Financial Economics}, 2(1):57--86, 2008.

\bibitem{englezos2009utility}
N.~Englezos and I.~Karatzas.
\newblock Utility maximization with habit formation: Dynamic programming and
  stochastic {PDE}s.
\newblock {\em Siam Journal on Control and Optimization}, 48(2):481--520, 2009.

\bibitem{he2013optimal}
L.~He and Z.~Liang.
\newblock Optimal investment strategy for the {DC} plan with the return of
  premiums clauses in a mean--variance framework.
\newblock {\em Insurance: Mathematics and Economics}, 53(3):643--649, 2013.

\bibitem{huang2018time}
Y.-J. Huang and A.~Nguyen-Huu.
\newblock Time-consistent stopping under decreasing impatience.
\newblock {\em Finance and Stochastics}, 22(1):69--95, 2018.

\bibitem{Karatzas2}
I.~Karatzas and S.~E. Shreve.
\newblock {\em Brownian Motion and Stochastic Calculus (Graduate Texts in
  Mathematics), 2nd edition}.
\newblock Springer, 1991.

\bibitem{karatzas1}
I.~Karatzas and S.~E. Shreve.
\newblock {\em Methods of Mathematical Finance (Stochastic Modelling and
  Applied Probability)}.
\newblock Springer, 1998.

\bibitem{kronborg2015inconsistent}
M.~T. Kronborg and M.~Steffensen.
\newblock Inconsistent investment and consumption problems.
\newblock {\em Applied Mathematics \& Optimization}, 71(3):473--515, 2015.

\bibitem{krusell}
K.~Krusell and A.~Smith.
\newblock Consumption-savings decisions with quasi-geometric discounting.
\newblock {\em Econometrica}, 71(1):365--375, 2003.

\bibitem{li2015time}
D.~Li, X.~Rong, and H.~Zhao.
\newblock Time-consistent reinsurance--investment strategy for an insurer and a
  reinsurer with mean--variance criterion under the {CEV} model.
\newblock {\em Journal of Computational and Applied Mathematics}, 283:142--162,
  2015.

\bibitem{li2013optimal}
Y.~Li and Z.~Li.
\newblock Optimal time-consistent investment and reinsurance strategies for
  mean--variance insurers with state dependent risk aversion.
\newblock {\em Insurance: Mathematics and Economics}, 53(1):86--97, 2013.

\bibitem{li2015timepartial}
Y.~Li, H.~Qiao, S.~Wang, and L.~Zhang.
\newblock Time-consistent investment strategy under partial information.
\newblock {\em Insurance: Mathematics and Economics}, 65:187--197, 2015.

\bibitem{pedersen2016optimal}
J.~L. Pedersen and G.~Peskir.
\newblock Optimal mean--variance selling strategies.
\newblock {\em Mathematics and Financial Economics}, 10(2):203--220, 2016.

\bibitem{pedersen2017optimal}
J.~L. Pedersen and G.~Peskir.
\newblock Optimal mean-variance portfolio selection.
\newblock {\em Mathematics and Financial Economics}, 11(2):137--160, 2017.

\bibitem{pollak1976habit}
R.~A. Pollak.
\newblock Habit formation and long-run utility functions.
\newblock {\em Journal of Economic Theory}, 13(2):272--297, 1976.

\bibitem{selten1965spieltheoretische}
R.~Selten.
\newblock Spieltheoretische behandlung eines oligopolmodells mit
  nachfragetr{\"a}gheit: Teil i: Bestimmung des dynamischen
  preisgleichgewichts.
\newblock {\em Journal of Institutional and Theoretical Economics}, (H.
  2):301--324, 1965.

\bibitem{selten1975reexamination}
R.~Selten.
\newblock Reexamination of the perfectness concept for equilibrium points in
  extensive games.
\newblock {\em International journal of game theory}, 4(1):25--55, 1975.

\bibitem{strotz}
R.~Strotz.
\newblock Myopia and inconsistency in dynamic utility maximization.
\newblock {\em The Review of Economic Studies}, 23(3):165--180, 1955.

\bibitem{weibull}
N.~Vieille and J.~W. Weibull.
\newblock Multiple solutions under quasi-exponential discounting.
\newblock {\em Economic Theory}, 39(3):513--526, 2009.

\bibitem{zeng2016robust}
Y.~Zeng, D.~Li, and A.~Gu.
\newblock Robust equilibrium reinsurance-investment strategy for a
  mean--variance insurer in a model with jumps.
\newblock {\em Insurance: Mathematics and Economics}, 66:138--152, 2016.

\end{thebibliography}
